\documentclass[11pt]{article}
 \usepackage{amssymb}
 \usepackage{epsfig}
 \usepackage{graphicx}
 \usepackage{pict2e}
 \newtheorem{Lemma}{Lemma}
 \newtheorem{Proposition}[Lemma]{Proposition}
 
 \newtheorem{Theorem}[Lemma]{Theorem}
 \newtheorem{Conjecture}[Lemma]{Conjecture}

  \newcommand{\BB}{\mbox{${\mathcal B}$}}
 \newcommand{\II}{\mbox{${\mathcal I}$}}
 \newcommand{\sfrac}[2]{{\textstyle\frac{#1}{#2}}}
 \newcommand{\proof}{{\bf Proof.\ }}
  
 \newcommand{\Reals}{{\mathbb{R}}}
  \newcommand{\Ints}{{\mathbb{Z}}}
  \newcommand{\bX}{{\mathbf X}}
    \newcommand{\bB}{{\mathbf B}}
   \newcommand{\qed}{\ \ \rule{1ex}{1ex}}
    \newcommand{\eps}{\varepsilon}
 
\begin{document}

\title{The Shape Theorem for Route-lengths in Connected Spatial Networks on Random
Points}
 \author{David J. Aldous\thanks{Department of Statistics,
 367 Evans Hall \#\  3860,
 U.C. Berkeley CA 94720;  aldous@stat.berkeley.edu;
  www.stat.berkeley.edu/users/aldous.  Aldous's research supported by
 N.S.F Grant DMS-0704159. }
}

 \maketitle
 
  \begin{abstract}
For a connected network on  Poisson points in the plane, consider the route-length
$D(r,\theta) $ between a point near the origin  and a point near polar
coordinates $(r,\theta)$, and suppose $E 
D(r,\theta) = O(r)$ as $r \to \infty$.
By analogy with the shape theorem for first-passage percolation, 
for a translation-invariant and ergodic network 
one expects $r^{-1} D(r, \theta)$ to converge as $r \to \infty$  to a constant
$\rho(\theta)$.
It turns out there are some subtleties in  making a precise formulation and a
proof.  We give one formulation and proof via a variant of the subadditive ergodic
theorem wherein random variables are sometimes  infinite.  
  \end{abstract}
 \vspace{0.1in}



 \vspace{0.4in}

 {\em Short title:}
Shape Theorem for Spatial Networks.

 \newpage
 \section{Introduction}
 \label{sec-INT}
 This paper is a technical part of a broader project investigating
connected random spatial networks, in particular networks built over a Poisson
process of points $\xi$ in the plane.  
 See \cite{me-spatial-1} for the least technical overview.    
 In any such network there is a (shortest-) route-length
$d(\xi,\xi^\prime)$ between each pair of points of the Poisson process, which by
connectivity is finite.  
 Under weak assumptions (see \cite{me-spatial-4} for a sufficient
condition) one expects the mean route-length to grow only linearly
with Euclidean distance.  
 Consider a (deliberately vague, for now) notion
 \begin{quote}
  $D(r,\theta) $ is the route-length between a point $\xi$ near the
origin 
and a point $\xi^\prime$ near polar coordinates $(r,\theta)$, 
 \end{quote}
 and suppose we know $E  D(r,\theta) = O(r)$ -- that is, suppose we have a linear upper
bound on mean route-length.  
 If the random network has translation-invariant and ergodic distribution,
then we intuitively expect that there should be a limit constant
 $\rho(\theta) = \lim_{r \to \infty} r^{-1} E  D(r,\theta)$ 
 and that in some sense renormalized random route-lengths should converge to the
limit constant:  
 $r^{-1}   D(r,\theta) \to \rho(\theta)$. 
This intuition arises in part from an analogy with the shape theorem  for
first-passage percolation \cite{MR905330,kesten-FPP} on the edges of the grid
$\Ints^2$.
In the usual such model the times $\tau(e)$ attached to edges $e$ are assumed
i.i.d., but the proof (based on the subadditive ergodic theorem) extends to the setting where the $\tau(e)$ are assumed only to be 
ergodic translation-invariant.   Studying route-lengths in random networks built
over Poisson point processes is perhaps the most natural continuum analog of studying
first-passage times in such lattice models. 
Two previously studied special continuum models,  superficially different, can be fitted into our general setup -- see section \ref{sec-analogs}.

 \subsection{Formulating a theorem}
 In the broader project we visualize a spatial network as having vertices and edges; in most
contexts, summary statistics such as 
 ``mean edge-length per unit area" are natural and important.  
 In the specific context of this paper, only the induced route-lengths
$d(\xi,\xi^\prime)$ are of interest, so 
we will dispense with other structure and work within the following set-up
throughout this paper.
 
 \noindent
 (A1) There is a Poisson process of points $\Xi = \{\xi\}$ of intensity
one, on $\Reals^2$. \\
 (A2) On each realization of $\Xi$ there are non-negative finite
``route-lengths" 
 $d(\xi,\xi^\prime) = d(\xi^\prime, \xi)$ which are assumed (only) to satisfy the triangle
inequality. \\
 (A3) The distribution of the whole structure $\{\xi; d(\xi, \xi^\prime) \}$ is
translation-invariant.
 That is, invariant w.r.t. the action of the group 
 $(T_{a,b}; \ (a,b) \in \Reals^2)$ where $T_{a,b}(x,y) = (x+a,y+b)$. 
 Moreover for each $(a,b) \neq (0,0)$ the action of $T_{a,b}$ is ergodic.
 
 Some specific examples are mentioned briefly in section \ref{sec-examples}, though our emphasis is on the generality of the assumptions.
 So it is worth mentioning what we are not assuming.  We are {\bf not} assuming 

 \noindent
 (B1) rotational invariance \\
 (B2) that  
 $d(\xi,\xi^\prime) \geq | \xi - \xi^\prime|$ (implicit in the underlying ``route-lengths" story)\\
 (B3)  any kind of  ``locality" for the route-lengths $d(\xi,\xi^\prime)$.
 
 \noindent   
In particular,  for a nearby pair $\xi,\xi^\prime$ the route-length 
$d(\xi,\xi^\prime)$ may depend on the entire configuration $\Xi$.
Finally, we often describe points in $\Ints^2$ by their radial coordinates. 
So $(r,\theta)$ denotes a point $z \in \Ints^2$; conversely, given $z \in \Ints^2$ we write $(r_z,\theta_z)$ for its radial coordinates.

 At first sight it looks easy to state and prove a theorem under assumptions (A1-A3) -- just find
a suitable formalization of the vague  notion $D(r,\theta)$ above, to which the subadditive
ergodic theorem can be applied.
But actually carrying this through seems surprisingly subtle.  
One attempt is to condition on points being planted at the origin and at
$(r,\theta)$:  the remaining points being still distributed as the Poisson point
process, one can define a conditioned network and then define $D(r,\theta)$ to be the
route-lengths between the planted points in the conditioned network.  However, for rather trivial reasons the desired result is simply not true
in this interpretation (see section \ref{sec-cex}). A second attempt is to 
interpret $D(r,\theta)$ as originally stated for the point {\em nearest} the origin
and the
point {\em nearest} $(r,\theta)$; this makes a precise definition 
 but it seems hard to work directly with this definition. 
A third attempt is to start by finding some feature to which one can apply the subadditive
ergodic theorem.  
 For instance, let $\xi_n$ be the leftmost point in the semi-infinite
strip 
 $ [n,\infty) \times [-1,1] $.  
 One can certainly apply the subadditive ergodic theorem to the array $(d(\xi_m,\xi_n))$ to conclude that
(under an integrability assumption) 
 $n^{-1} d(\xi_0,\xi_n)$ converges a.s. to a constant limit $\rho(0)$.  
At first sight this approach seems to resolve the whole issue.  
But the precise hypotheses and conclusions are tied to the particular feature
initially chosen, and it seems technically hard to reconcile the results from
different choices.

We adopt a fourth approach, aimed at a more natural type of conclusion. Write $A,B$
for bounded subsets of $\Reals^2$ and for $z \in \Reals^2$  write $z + B =
\{z+z^\prime: \ z^\prime \in B\}$.  
To motivate the precise definitions (\ref{def-S},\ref{LSP}) below, 
consider $\xi \in A$ and $\xi^\prime \in z+B$ with $r_z$ large; 
the route-length $d(\xi,\xi^\prime)$ provides one interpretation of our initial vague notion $D(r_z,\theta_z)$, which we want to prove is 
approximately the (deterministic) length 
$r_z \rho(\theta_z)$.  
To avoid conditioning on existence of points in sets, we sum: 
$\sum_{\xi \in A} \sum_{\xi^\prime \in z+B} d(\xi,\xi^\prime)$ should be  $N(A) N(z+B) r_z \rho(\theta_z) \pm o(r_z)$ where $N(\cdot)$ is the counting process of $\Xi$.
But we can avoid writing $N(\cdot)$ by rewriting the approximation as 
$\sum_{\xi \in A} \sum_{\xi^\prime \in z+B} |d(\xi,\xi^\prime) - r_z \rho(\theta_z)| = o(r_z)$.  This prompts the following definitions. 

 For $c \geq 0$ define a
random variable 
\begin{equation}
S(A,B; c) := \sum_{\xi \in A} \sum_{\xi^\prime \in B}
|d(\xi,\xi^\prime) - c|  . 
\label{def-S}
\end{equation}
Say the random network has the 
{\em $L^1$ shape property} if there exist constants $\rho(\theta)$  such that, for
all bounded $A,B$,
\begin{equation}
r_z^{-1} E S(A, z + B; r_z \rho(\theta_z)) \to 0 \mbox{ as } r_z \to \infty  .
\label{LSP}
\end{equation}
Because $S$ is an additive set function, it is enough to prove (\ref{LSP}) when $A$
and $B$ are sufficiently large (or sufficiently small) squares centered at the origin;
and in the latter case we see how this notion provides a formalization of the idea behind $D(r,\theta)$. 

Having decided on the conclusion we seek, what hypotheses do we need?   
Obviously it is necessary that the corresponding linear upper bound holds: for all
bounded $A,B$,
\begin{equation}
E \sum_{\xi \in A} \sum_{\xi^\prime \in z+B}
d(\xi,\xi^\prime)
 = O(r_z) \mbox{ as } r_z \to \infty  . \label{LSP-ub}
\end{equation}
We conjecture that (\ref{LSP-ub}) is sufficient (see section \ref{sec-as-shape} for precise statement). 
However in this paper we work under the analogous, but stronger,  $L^2$ assumption:
for all bounded $A,B$,
\begin{equation}
\sup_z \frac{
E \sum_{\xi \in A} \sum_{\xi^\prime \in z+B}
d^2(\xi,\xi^\prime)}
{\max(1,r^2_z)} < \infty .
\label{Hyp}
\end{equation}
Again, it is enough to verify this when $A$
and $B$ are sufficiently  small squares centered at the origin.

\begin{Theorem}
\label{T1}
Under the standing assumptions (A1 - A3), if hypothesis (\ref{Hyp}) holds 
then the $L^1$ shape property (\ref{LSP}) holds. 
Moreover
\begin{equation}
\sup_{\theta_2 \neq \theta_1} \frac{|\rho(\theta_2) - \rho(\theta_1)|}{|\theta_2-\theta_1|} < \infty .
\label{Lip}
\end{equation}
\end{Theorem}
This is proved in section \ref{sec-proof}, though the main work is delegated to a new 
``subadditive ergodic theorem with missing values", Proposition \ref{P-sub}, stated and proved in section \ref{sec-subadd}.
A conjectured stronger ``a.s. shape theorem" conclusion is discussed briefly in section \ref{sec-as-shape}.
Obviously, if we add the assumption of rotational invariance then $\rho(\theta)$ is constant.

\subsection{Examples using route-lengths}
 \label{sec-examples}
In all these examples, the $d(\xi,\xi^\prime)$ are minimal 
route-lengths within given networks.

\paragraph{Proximity graphs \cite{jaromczyk,me-spatial-1}.}
This family of graphs (our main example) is defined by: 
\begin{quote}
$(\xi,\xi^\prime)$ is an edge iff 
the set $A(\xi,\xi^\prime)$ contains no other point of the Poisson process 
\end{quote}
for different choices of $A(\xi,\xi^\prime)$, a fundamental choice 
(giving the {\em relative neighborhood graph}) being
\begin{quote}
 $A(\xi,\xi^\prime)$ is the intersection of 
the disc with center $\xi$ and radius $|\xi^\prime - \xi|$ 
and the disc with center $\xi^\prime$ and radius $|\xi^\prime - \xi|$.
\end{quote} 
Other graphs in the family use subsets of this $A$ and hence are supergraphs of the
relative neighborhood graph and hence can only have smaller route-lengths.
The purpose of the companion paper \cite{me-spatial-4} is to 
give a general property that implies our present  condition (\ref{Hyp}), 
and to verify this property for the relative neighborhood graph.  
It follows that (\ref{Hyp})  holds for the relative neighborhood graph 
(and hence for every proximity graph) on a Poisson point process.
Because edges are defined by a deterministic rule, proximity graphs 
inherit the stationary ergodicity property (A3) from the trivial tail $\sigma$-field
property of the Poisson process.
So our Theorem \ref{T1} applies, and by
rotational invariance $\rho(\theta)$ is a constant $\rho$, depending on the model. 
Monte Carlo estimates of $\rho$ (around 1.4 for the relative neighborhood graph) can be seen in 
\cite{me-spatial-1} but we do not know any explicit rigorous upper bound.  

The general condition in \cite{me-spatial-4} might be applicable to other models, but the examples below can be 
handled more directly.
 
\paragraph{Lattice-based networks.}
One can start with (for instance) the square grid lattice as a network, and simply 
connect each Poisson point to the nearest grid point.  
One can see directly that this random network satisfies the shape property with
$\rho(\theta) = |\cos \theta| +  | \sin \theta |$.  
This conclusion remains true if we make the network be translation-invariant and ergodic by replacing 
the deterministically-spaced grid lines by randomly-spaced ones.

\paragraph{Asymptotically efficient networks.}
It is not surprising that there are networks which are 
``optimal" in the sense $\rho(\cdot) \equiv 1$. 
It is at first sight surprising that one can find such networks whose
length-per-unit-area 
is arbitrarily close to the minimum possible (over all connected 
networks -- attained by the Steiner tree) length-per-unit-area.
But this can be achieved by the simple device of superimposing, over the Steiner tree, 
a sparse Poisson line process.
This construction is studied in detail in \cite{me116}.

\paragraph{The Hammersley network.} 
This network, introduced in \cite{me-spatial-1}, has the remarkable
property that at each point $\xi$ there are exactly $4$ edges, one in each of the
four quadrant 
directions (i.e. between East and North, etc).  This network has not been studied
carefully, 
but it is plausible one can use known properties of the underlying Hammersley
process to prove directly
that the shape property holds with 
$\rho(\theta) = \rho_0 (|\cos (\theta - \pi/4)| +  | \sin (\theta - \pi/4)|)$ 
for some constant $\rho_0$.

\subsection{Other examples}
\label{sec-analogs}
Suppose we remove the ``satisfy the triangle inequality" requirement from (A2), 
to get instead
\begin{quote}
 (A2*) On each realization of $\Xi$ there are ``costs" 
 $0 < c(\xi,\xi^*) = c(\xi^*, \xi) \le \infty$.
\end{quote}
One can now define $d(\xi,\xi^\prime)$ as the cost of the minimum-cost path from $\xi$ to $\xi^\prime$, and this makes $d$ satisfy the triangle inequality.  So, provided $d$ is always finite, (A2) holds.
If the other hypotheses of Theorem \ref{T1} hold for $d$, then the conclusion of Theorem \ref{T1} gives the $L^1$
shape property
for $d$.  
The following two particular cases have been studied previously by direct methods
which establish 
the a.s. shape theorem; our Theorem \ref{T1} applies (assuming second moments in (a)) to
give the $L^1$ 
shape theorem.

% Our model -- route-lengths in translation-invariant random networks built
% over Poisson point processes -- differs from two previously proposed 
% ``continuum analogs of first-passage percolation", also built over Poisson processes.

 \noindent
(a) Take the Delaunay triangulation on the Poisson points, and then take
$c(\xi,\xi^\prime)$ to be i.i.d. with finite mean on the 
edges of the triangulation (and $= \infty$ elsewhere): \cite{MR1166620}.

\noindent
(b) Fix  $\alpha > 1$  and set $c(\xi,\xi^\prime) = |\xi - \xi^\prime|^\alpha$: 
\cite{MR1452554}.

\section{Reducing the proof of Theorem \ref{T1} to a subadditivity result}
\label{sec-proof}
Write $A, B$ for  bounded subsets of  $\Reals^2$ with non-zero area. 
Write $N(A)$ for the number of points of the Poisson point process $\Xi$ in $A$ 
and write $G(A)$ for the ``good" event $\{N(A) \geq 1\}$.
On $G(A)$ let $\xi_A$ be a uniform random point of 
$\Xi \cap A$.  
Note that hypothesis (\ref{Hyp}) implies
\begin{equation}
E[ d(\xi_{z_1+A},\xi_{z_2 + B}) 1_{G(z_1+A) \cap G(z_2+B)}]  \leq \kappa(A,B) \ {\max(1,|z_1 - z_2|)} 
\label{eq-kappa}
\end{equation}
where $\kappa(A,B) < \infty$ depends only on $A,B$.

Fix $\theta \in [0,2\pi)$ and fix a bounded subset $A \subset \Reals^2$ of non-zero area. 
Write $z_n$ for the point with radial coordinates $(nr_0,\theta)$, 
where $r_0$ is sufficiently large that the sets $z_n + A$ are disjoint.
So $ G(z_n + A)$ is the event $\{N(z_n + A) \geq 1\}$, and on $G(z_n + A)$ let $\xi_{z_n+A}$ be a uniform random point of 
$\Xi \cap (z_n + A)$.

Consider the array of random variables
\begin{eqnarray}
X_{mn} &=& d(\xi_{z_m+A},\xi_{z_n+A}) \mbox{ on } G(z_m+A) \cap G(z_n+A) \label{star} \\
              &=& \infty \mbox{ otherwise } . \nonumber
\end{eqnarray}
Note (A2) implies the triangle inequality 
\begin{equation}
\mbox{for $\ell<m<n$,  $X_{\ell n} \leq X_{\ell m} + X_{mn} \mbox{  on } 1_{G(z_\ell +A) \cap
G(z_m +A) \cap G(z_n +A)}$ .} 
\label{tria}
\end{equation}
Proposition \ref{P-sub}, stated and proved in section \ref{sec-subadd}, is tailored to this setting.
Specifically, hypothesis (i) is (\ref{star}), (ii) is (\ref{tria}), (iii) follows from Poisson independence, (iv) from (A3) and (v) from (\ref{Hyp}). 
 The conclusion of Proposition \ref{P-sub} is that
there exists a constant $0 \leq \rho(\theta) < \infty$ such that 
\begin{equation}
   E\left[ \left| \frac{d(\xi_A,\xi_{z_n+A})}{nr_0} \ -\rho(\theta) \right| 1_{G(A) \cap G(z_n+A)} \right] \to 0 .
\label{rho-1}
\end{equation}
This is the main ingredient of the proof; the argument below continues with the details of converting (\ref{rho-1}) into the stated conclusion of
 Theorem \ref{T1}. 
The typography in (\ref{rho-1}) is potentially confusing; note we are multiplying an absolute value by an indicator, not taking a conditional expectation.

{\em A priori} the limit constant $\rho(\theta)$ in (\ref{rho-1}) might depend on $A$ and on $r_0$.  
We first show it does not depend on $r_0$; more precisely, we will show
\begin{equation}
   E\left[ \left| \frac{d(\xi_A,\xi_{(r,\theta)+A})}{r} \ -\rho(\theta) \right| 1_{G(A) \cap G((r,\theta)+A)} \right] \to 0  \mbox{ as } r \to \infty .
\label{rho-2}
\end{equation}
Let us give the argument for (\ref{rho-2}) in some detail, intending to omit similar details later.  
Write $r = (n_1 + n_2)r_0 + \gamma$ for some $0 \leq \gamma < r_0$. 
Let $r \to \infty$ while choosing $n_1 = n_1(r) \to \infty$ and $n_2 = n_2(r) \to \infty$.  
By (\ref{rho-1}) and translation-invariance
\[
   E\left[ \left| \frac{d(\xi_A,\xi_{(n_1r_0,\theta)+A}) - n_1r_0 \rho(\theta)}{r}  \right| 1_{G(A) \cap G((n_1r_0,\theta)+A)} \right] \to 0   
   \]
   \[
    E\left[ \left| \frac{d(\xi_{(n_1r_0 + \gamma,\theta)+A},\xi_{(r,\theta)+A}) - n_2r_0 \rho(\theta)}{r} \ \right| 1_{G((n_1r_0+\gamma,\theta)+A
) \cap G((r,\theta)+A)} \right] \to 0  .
  \] 
  Combining these with (\ref{eq-kappa}) applied to $\xi_{(n_1r_0,\theta)+A}$ and $\xi_{(n_1r_0 + \gamma,\theta)+A}$, and using the 
  triangle inequality for 
  $d(\cdot,\cdot)$, 
  \[
     E\left[ \left| \frac{d(\xi_A,\xi_{(r,\theta)+A})}{r} \ -\rho(\theta) \right| 1_{G(A) \cap  G((r,\theta)+A) \cap  G((n_1r_0,\theta)+A)  \cap 
G((n_1r_0+\gamma,\theta)+A)   } \right] \to 0  \mbox{ as } r \to \infty .
   \]
  This expression differs  from (\ref{rho-2}) only by the inclusion of the restriction to $G((n_1r_0,\theta)+A)  \cap G((n_1r_0+\gamma,\theta)+A)
$, an event which has probability at least $1 - p$ for some $p = p(r_0,A) < 1$.  Given $k$ we can choose (for large $r$) $k$ different values of 
$n_1$ such that the $k$ corresponding events are independent because the underlying sets are disjoint; it follows that 
    \[
     E\left[ \left| \frac{d(\xi_A,\xi_{(r,\theta)+A})}{r} \ -\rho(\theta) \right| 1_{G(A) \cap  G((r,\theta)+A) \cap 
      H(r,k) } \right] \to 0  \mbox{ as } r \to \infty 
   \]
   for certain events $H(r,k)$ such that $P(H(r,k)) \geq 1 - p^k$ for large $r$.
   Letting $k \to \infty$ and appealing to the $L^2$ bound (\ref{Hyp}) establishes (\ref{rho-2}).
   
   Now consider two subsets $A \subset A^\prime$.  
   Could the two constants $\rho(\theta)$ and $\rho^\prime(\theta)$ in  (\ref{rho-2}) be different?  
   When we make independent  choices of random points $\xi_{z+A}$ and $\xi_{z + A^\prime}$ there is some fixed probability 
   $p(A,A^\prime) > 0$ that the two random points are the same, and it easily follows that the limit constants must be equal.
   That is, $\rho(\theta)$ does not depend on choice of $A$.
   
   Next we prove the Lipschitz property (\ref{Lip}).  
   Fix $\theta_1$ and $\theta_2$.
   The triangle inequality and (\ref{eq-kappa}) give 
   \[ E[ |d(\xi_A,\xi_{(r,\theta_1)+A}) - d(\xi_A,\xi_{(r,\theta_2)+A})| 1_{G(A) \cap G((r,\theta_1)+A) \cap G((r,\theta_2)+A)}]  \leq \]
    \[ E[ d(\xi_{(r,\theta_1)+A} ,\xi_{(r,\theta_2)+A})| 1_{G(A) \cap G((r,\theta_1)+A) \cap G((r,\theta_2)+A)}]  \leq 
    \kappa(A,A) \min (1, r |\theta_2-\theta_1|) .\]
    Applying (\ref{rho-2}),
    \[ |\rho(\theta_2) - \rho(\theta_1)| P^3(G(A)) \leq  \kappa(A,A) \ |\theta_2-\theta_1| \]
  and now any choice of $A$  establishes (\ref{Lip}).

   Next we want to prove the analog of (\ref{rho-2}) where the angle is not fixed.  
   That is, for $z = (r_z,\theta_z)$ with $r_z \to \infty$, we claim
   \begin{equation}
   E\left[ \left| \frac{d(\xi_A,\xi_{z+A})}{r_z} \ -\rho(\theta_z) \right| 1_{G(A) \cap G(z+A)} \right ] \to 0  \mbox{ as } r_z \to \infty .
\label{rho-3}
\end{equation}
By compactness and continuity of $\rho(\cdot)$ we can reduce to the case $(r_n,\theta_n)$ where $\theta_n \to \theta$, and it is enough to prove 
   \begin{equation}
   E\left[ \left| \frac{d(\xi_A,\xi_{(r_n,\theta_n)+A})}{r_n} \ -\rho(\theta) \right| 1_{G(A) \cap G((r_n\theta_n)+A)} \right ] \to 0  \mbox{ as 
} r_n \to \infty .
\label{rho-4}
\end{equation}
Here we repeat the format of the argument for (\ref{rho-2}).
Take $r^*_n = (1+o(1))r_n$.   Apply the triangle inequality to $\xi_A, \xi_{(r^*_n,\theta)+A}, \xi_{(r_n,\theta_n)+A}$, apply the fixed-$\theta$
 result (\ref{rho-2}) to 
the first distance and apply (\ref{eq-kappa}) to the second distance; 
we deduce the analog of (\ref{rho-4}) with the extra term $1_{G((r^*_n,\theta)+A)}$. 
But this is true for each of  multiple possible choices for $r^*_n$, so we can deduce (\ref{rho-4}) and thence (\ref{rho-3}).  

To complete the proof we need to convert (\ref{rho-3}) into an assertion involving the sums $S(A, z + B; r_z \rho(\theta_z))$ appearing in (\ref{LSP}).  
Let us state the underlying logical structure carefully; note there is no assumption that the $(Y^{(n)}_{ij})$ are independent of  $(N^{(n)}_1, N^{(n)}_2)$.
\begin{Lemma}
\label{L1}
Fix $\lambda_1, \lambda_2$. 
For each $n$ let 
$(Y^{(n)}_{ij}, 1 \leq i \leq N^{(n)}_1, 1 \leq j \leq N^{(n)}_2)$ be an array of nonnegative random variables, and suppose that $N^{(n)}_1$ and 
$N^{(n)}_2$ are independent with Poisson($\lambda_1$) (resp. $\lambda_2$) distributions.  
On the event $\{ N^{(n)}_1 \geq 1, N^{(n)}_2 \geq 1\}$, and conditional on the entire collection $(Y^{(n)}_{ij}, 1 \leq i \leq N^{(n)}_1, 1 \leq 
j \leq N^{(n)}_2)$, take 
$(U^{(n)}_1,U^{(n)}_2)$ to be independent with Uniform$[1,2,\ldots, N^{(n)}_1]$ and Uniform$[1,2,\ldots, N^{(n)}_2]$ distributions.  
Suppose
\[ Y^{(n)}_{U^{(n)}_1,U^{(n)}_2} 
\ 1_{(N^{(n)}_1 \geq 1, N^{(n)}_2 \geq 1)} 
\to 0 \mbox{ in probability as } n \to \infty . \]
Then 
\[ \sum_{i=1}^{N^{(n)}_1} \sum_{j=1}^{N^{(n)}_y} \ Y^{(n)}_{ij}   \to 0 \mbox{ in probability as } n \to \infty . \]
\end{Lemma}
\proof
It is enough to prove the conclusion restricted to $\{ 1 \leq N^{(n)}_1 \leq L, 1 \leq N^{(n)}_2 \leq L\}$ for fixed $L$.
But with this restriction, the hypothesis implies 
$\max_{ij} Y^{(n)}_{ij} \to 0$ in probability, which in turn implies the conclusion. 
\qed

\vspace{0.2in}
Now consider the setting of (\ref{rho-3}).  As $r_z \to \infty$ 
the array 
\[  \left(\left| \frac{d(\xi,\xi^\prime)}{r_z} \ -\rho(\theta_z) \right| , \ \xi \in A, \ \xi^\prime \in z+A \right) \] 
satisfies the assumptions of Lemma \ref{L1}, and the conclusion is 
\[ 
r_z^{-1} \  S(A, z + A; r_z \rho(\theta_z)) \to 0 \mbox{ in probability}. \] 
The $L^2$ bound (\ref{Hyp}) extends this to 
\[r_z^{-1} \ E  S(A, z + A; r_z \rho(\theta_z)) \to 0 \] 
which is enough to establish the $L^1$ shape property.
  
\subsection{The conjectured a.s. shape theorem}
\label{sec-as-shape}
Implicit in the underlying picture of route-lengths in spatial networks is that route-lengths are at least as big as Euclidean distance: 
\begin{equation}
d(\xi,\xi^\prime) \geq |\xi - \xi^\prime| .\label{Euclid} 
\end{equation}
This was not assumed for Theorem \ref{T1}; assuming it here, we see $\rho(\theta) \geq 1$. 
Using  ``triangle inequality" arguments as in the previous section, it is easy to check that 
\[ \bB :=  \{z = (r,\theta): \ r \leq 1/\rho(\theta) \} \]
defines a {\em convex} subset of the unit disc.  
A natural informal statement of a shape theorem is that, if we plant one Poisson point $\xi_0$ at the origin, then for large $\ell$ the set of 
points at route-length at most $\ell$ from $\xi_0$ is approximately the set of points within $\ell \bB$.  
So one can formalize the {\em a.s. shape property} as follows, in the context of a planted point $\xi_0$ at the origin.
For each $\eps > 0$ there exists random $L(\eps) < \infty$ such that for all $\ell > L(\eps)$
\[
\Xi \cap (1 - \eps) \ell \bB \subseteq \{\xi \in \Xi: \ d(\xi,\xi_0) \leq \ell \} \subseteq (1+\eps) \ell \bB . \] 
This has been proved by direct methods in the two special models of section \ref{sec-analogs}.
Because our Theorem \ref{T1} conclusion involves $L^1$ convergence instead of a.s. convergence, it implies only a somewhat  weaker result; and also our `$`L^2$ bounded" assumption is stronger than seems necessary.  
In other words, the natural conjecture suggested by the analogy with the  
shape theorem for first-passage percolation is as follows.
\begin{Conjecture}
Under the standing assumptions (A1 - A3), and (\ref{LSP-ub}) and (\ref{Euclid}), there exists a convex set $\bB$ such that the a.s. shape property holds.
\end{Conjecture}
And though we work throughout with an underlying Poison point process, such a result might be expected to hold for any ergodic translation-invariant point process.

\subsection{A minor counter-example}
\label{sec-cex}
Take a model to which Theorem \ref{T1} applies with $\rho(\theta) > 1$.  
Choose $r_k \uparrow \infty$ fast and $\delta_k \downarrow 0$. 
Modify the model by putting a straight line link between each pair of points whose distance apart is in 
$\cup_k [r_k,r_k+\delta_k]$.   
By making the $\delta_k \downarrow 0$ sufficiently fast and appealing to the ``minimum cost path" device of section \ref{sec-analogs},
 the hypotheses and conclusion of Theorem \ref{T1} remain true with the same limit $\rho(\theta)$ -- the extra links make no difference to route-length between typical pairs.
 But if we had attempted to formulate the theorem using 
 ``$D(r,\theta) = $ distance between two points at distance $r$ apart" then we would not get a $r \to \infty$ limit, because 
 of the  exceptional $r \in \cup_k [r_k,r_k+\delta_k]$.

 \section{A subadditive ergodic theorem with missing values}
 \label{sec-subadd}
 We develop a variation of Kingman's subadditive ergodic theorem (see e.g.
\cite{MR0438477}) in which the random variables $(X_{ij}, 0 \leq i < j < \infty)$
are sometimes undefined 
 (in which case we will set the value to $\infty$, though that isn't quite
the natural interpretation in our application).
 Consider a sequence $(G_i, 0 \leq i < \infty)$ of ``good" events, and
write 
 $I_i$ for the indicator $1_{G_i}$ and write $I_{ij} = I_iI_j = 1_{G_i
\cap G_j}$.  
 Our assumptions are

 \noindent
 (i) $0 \leq X_{ij} \leq \infty; \quad X_{ij} < \infty \mbox{ on } G_i \cap G_j
$ .\\
 (ii) For $i<j<k$,  $X_{ik} \leq X_{ij} + X_{jk} \mbox{  on } G_i \cap
G_j \cap G_k$ .\\
 (iii) The process $(I_i, 0 \leq i < \infty)$ is independent
Bernoulli($\delta$) for fixed $0<\delta<1$.\\
 (iv) Setting $\bX_i = (I_i, X_{i,i+k}, 1 \leq k< \infty)$, the process 
 $(\bX_i, 0 \leq i < \infty)$ is stationary and ergodic.\\
 (v) $ \sup_{n \geq 1} n^{-2} E [X^2_{0n} I_{0n} ]< \infty$.
 
 Mostly these are the obvious analogs of the usual assumptions
\cite{MR0438477}.
 Note that in the usual setting we have a trivial implication
 \[ \mbox{ if } EX_{01} < \infty \mbox{ then  } \sup_n n^{-1} EX_{0n} \le EX_{01} <
\infty \]
 whereas in our setting the implication:  
  \[ \mbox{ if } E[X_{01}I_{01} ] < \infty \mbox{ then  } \sup_n n^{-1}
E[X_{0n}  I_{0n}] < \infty \]
  is not trivial (we don't know if it is true). 
  The latter would be the natural hypothesis in our setting, but to make
our straightforward proof technique work we make the stronger $L^2$ assumption (v).  
  Also to keep matters simple, we assume ergodicity and seek only $L^1$
convergence.
  \begin{Proposition}
  \label{P-sub}
  Assume (i)-(v).
  Then there exists a constant $0 \leq c < \infty$ such that 
  $E\big[ | \frac{X_{0n}}{n} \ -c | I_{0n} \big] \to 0 $.
  \end{Proposition}
  \proof 
  We compare the given process with another process in which one is
allowed to use the ``bad" indices, but with high penalty.  
  Fix large $K$.
  Define
   \begin{eqnarray*}
   \widetilde{X}_{ij} &=& X_{ij} \mbox{ on } G_i \cap G_j\\
     \widetilde{X}_{i,i+1} &=& K \mbox{ on the complement } (G_i \cap G_{i+1})^c
\mbox{ of } G_i \cap G_{i+1} \\
    \widetilde{X}_{ij} && \mbox{ undefined, otherwise. }   
  \end{eqnarray*}
  Now define a process
  \[ Y_{ij} = \min \left (  \widetilde{X}_{i_0,i_1} + 
\widetilde{X}_{i_1,i_2}  + \ldots +  \widetilde{X}_{i_{m-1},i_m} 
\right) ,\] 
  the minimum over $i = i_0<i_1< i_2 < \ldots <i_m = j$ such that each
$\widetilde{X}_{i_{u-1},i_u} $ is defined.  
  Observe that $Y_{ij} $ is always defined and finite, and is subadditive.
 Also, because 
  $Y_{01} = X_{01}I_{01} + K(1 -  I_{01} )$ we have $E Y_{01} < \infty$. So we
can apply Kingman's subadditive ergodic theorem to deduce there
exists a constant 
  $0 \leq c^{(K)} < \infty$ such that 
  \begin{equation}
   E  \left| \frac{Y^{(K)}_{0n}}{n} \ -c^{(K)} \right| \to 0  \mbox{ as } N \to
\infty   
   \label{K1}
   \end{equation}
  where we now write $Y^{(K)}$ to emphasize dependence on $K$.
  
  Note that the $L^1$ convergence in (\ref{K1}) implies
  \[ \delta^2  c^{(K)} = \lim_{n \to \infty} E[ c^{(K)} I_{0n} ]  =
\lim_{n \to \infty} E[ \sfrac{Y^{(K)}_{0n}}{n}     I_{0n} ]  . \]
  Now by assumption (v) we have
  \[ B_1 := \sup_n n^{-1} E [X_{0n} I_{0n}] < \infty \]
  and from the  definition 
  \begin{equation}
   Y^{(K)}_{0n} \leq X_{0n} \mbox{ on } G_0 \cap G_n. 
   \label{YKX}
   \end{equation}
  So $\delta^2  c^{(K)} \leq B_1$. 
  Also from the definition we see that $Y^{(K)}_{0n}$ is non-decreasing in
$K$.
  Hence so is $c^{(K)}$, and so  we can define the limit 
  \[ c:= \lim_{K \to \infty} c^{(K)}  \]
  for which $0 \leq c < \infty$.
  
  The essential issue is to bound the difference in (\ref{YKX}).   Fix $K$
and $n$.
  By definition there is some path $0 = i_0<i_1< i_2 < \ldots <i_m = n$
for which 
  \begin{equation}
    Y^{(K)}_{0n} =  \widetilde{X}_{i_0,i_1} +  \widetilde{X}_{i_1,i_2}  +
\ldots +  \widetilde{X}_{i_{m-1},i_m}  . 
    \label{sum}
    \end{equation}
  Some steps on the path are of the form (for some $i$)
  \[
  (i_u,i_{u+1}) = (i,i+1) \mbox{ and } G_i \cap G_{i+1} \mbox{ fails, so }
\widetilde{X}_{i_{u},i_{u+1}}  = K . \]
  Write $\II$ for the random set of $i$ for which this occurs (for some
$i_u$), and note for later use 
  \begin{equation}
   Y^{(K)}_{0n} \geq K |\II |  . \label{II-def}
  \end{equation}
  
  Now consider  a maximal run $[b,b^\prime]$ of bad events; that is, $G_i$
occurs for $i = b-1$ and $i = b^\prime + 1$ but not for $b \leq i \leq b^\prime$.  
  Then the path in (\ref{sum})  must either use all the edges 
  $(b-1,b), \ (b,b+1), \ \ldots , (b^\prime -1,b^\prime), \
(b^\prime,b^\prime + 1)$
  or none of them.  
  If it uses all of them, replace the path segment $b-1\to b \to \ldots 
\to b^\prime \to b^\prime + 1$ by the single edge $b-1 \to  b^\prime + 1$, that is
replace the part 
  $\widetilde{X}_{b-1,b}  + \ldots +  \widetilde{X}_{b^\prime ,b^\prime
+1}$ of the sum (\ref{sum}) by $X_{b-1,b^\prime +1}$.  
  Make this replacement for each bad run touched by the path, and assume
we are on $G_0 \cap G_n$ so there are no endpoint issues.  This converts (\ref{sum}) into
a new sum 
   \begin{equation}
    Z^{(K)}_{0n} =  X_{j_0,j_1} +  X_{j_1,j_2}  + \ldots + 
X_{j_{q-1},j_q}   
    \label{sum2}
    \end{equation}
    where all the steps are between good indices, and so by the
subadditivity assumption (ii) we have 
    \[ 
    X_{0n} \leq Z^{(K)}_{0n} \mbox{ on } G_0 \cap G_n . \]
    The net effect of this conversion can be written precisely as  \[
Z^{(K)}_{0n} -  Y^{(K)}_{0n} = \sum_{(b,b^\prime) \in \BB}
X_{b-1,b^\prime +1} - K |\II|    \mbox{ on } G_0 \cap G_n \]
    where $\BB$ is the set of bad runs touched by the path. 
    So  we can bound the difference in (\ref{YKX}) rather crudely as 
       \begin{equation} 
    E \left[( X_{0n}  - Y^{(K)}_{0n} ) I_{0n} \right] \leq E \left[   I_{0n}   
\sum_{(b,b^\prime) \in \BB} X_{b-1,b^\prime +1} \right] 
    . \label{sum3}
    \end{equation}
  Recall hypothesis (v): 
   \[ B_2 := \sup_m m^{-2} E [X^2_{0m} I_{0m}] < \infty  .\]
  For $0 \leq i,j \leq n, \ j-i \geq 2$ write 
  $\Lambda_{ij}$ for the event $\{(i+1,j-1) \in \BB\}$ 
  and note $\Lambda_{ij} \subset I_{ij}$.  
  Using the Cauchy-Schwarz inequality.
   \begin{eqnarray}
   E [1_{ \Lambda_{ij}} X_{ij} I_{0n} ] 
   &\leq& \sqrt{P(  \Lambda_{ij} \cap I_{0n}) } \ \sqrt{ E
[X^2_{ij}I_{ij}]} \nonumber\\
     &\leq& \sqrt{P(  \Lambda_{ij} \cap I_{0n}) } \ (j-i) B_2^{1/2}
\label{sum4}
   \end{eqnarray} 
  the second inequality by the stationarity assumption (iv). 
  
  Set $p_{ij} = P(  \Lambda_{ij} \cap I_{0n})$ and note
  \[ p_{ij} \leq P(  \Lambda_{ij} ) \leq P(G^c_{i+1} \cap \ldots \cap
G^c_{j-1}) = (1-\delta)^{j-i-1} \]
  by assumption (iii).  
  Because the sum in (\ref{sum3}) can be written as 
  $\sum_{ij} 1_{\Lambda_{ij}} X_{ij} $, we can combine (\ref{sum3}) and
(\ref{sum4}) to get 
  \begin{equation} 
 E [( X_{0n}  - Y^{(K)}_{0n} ) I_{0n} ] \leq B^{1/2}_2  \sum_{i=0}^{n-2}
\sum_{j=i+2}^n (j-i) p_{ij}^{1/2} .
 \label{qwe-1}
     \end{equation}
 Now consider the double sum above with $p_{ij}^{1/2}$ replaced by
$p_{ij}$.  That is, consider 
 \begin{eqnarray}
  \sum_{i=0}^{n-2} \sum_{j=i+2}^n (j-i) p_{ij}
  &=& E I_{0n}   \sum_{i=0}^{n-2} \sum_{j=i+2}^n (j-i) 1_{ \Lambda_{ij}}
\label{qwe-0}\\
  &=& E I_{0n} |\II | \mbox{ for the random set $\II$ in
(\ref{II-def})}\nonumber\\\
  &\leq& K^{-1}  E I_{0n} Y^{(K)}_{0n} \mbox{ by
(\ref{II-def})}\nonumber\\\
  &\leq& K^{-1}  E I_{0n} X_{0n} \mbox{ by (\ref{YKX})}\nonumber\\\ &\leq&
B_1 n/K . \label{qwe}
 \end{eqnarray} 
 Now an elementary inequality (stated and proved as Lemma \ref{L2} below)
bounds the right side of (\ref{qwe-1}) in terms of the left side of  (\ref{qwe-0}).
 Combining this inequality with (\ref{qwe-1},\ref{qwe}) gives: for any $J
\geq 2$, 
 \[
 B^{1/2}_2  n^{-1}  E [( X_{0n}  - Y^{(K)}_{0n} ) I_{0n} ] \leq 
  \sum_{j=J+1}^\infty j (1 - \delta)^{(j-1)/2} 
  \ + \ J  \sqrt{B_1/K}
    . \]
Taking $J = K^{1/2}$ we see 
\[ \lim_{K \to \infty} \sup_n n^{-1}  E [( X_{0n}  - Y^{(K)}_{0n} ) I_{0n} ] = 0 \]
which by the $L^1$ convergence  (\ref{K1})  and the inequality 
(\ref{YKX})   implies 
\[ \lim_{K \to \infty} \limsup_n  E [ |n^{-1} X_{0n}  - c^{(K)}| I_{0n} ] = 0 \] 
establishing  Proposition \ref{P-sub} for $c= \lim_{K \to \infty} c^{(K)} $ which was
previously shown to be finite.

 \begin{Lemma}
 \label{L2}
 Let $0<\eta < 1$ and let 
 $(p_{ij}, 0 \leq i,j \leq n, \ j-i \geq 2)$ 
 be constants such that
 $0 \leq p_{ij} \leq \eta^{j-i-1}$.  Then for any $J \geq 2$
 \[ n^{-1} \sum_{i=0}^{n-2} \sum_{j=i+2}^n (j-i) p_{ij}^{1/2} \leq 
\sum_{j=J+1}^\infty j \eta^{(j-1)/2} \ + \ 
J  n^{-1/2} \sqrt{ \sum_{i=0}^{n-2} \sum_{j=i+2}^n (j-i) p_{ij} }. \]
 \end{Lemma}
 \proof
 Fix $i$ and set $q_j = p_{i,i+j}$ for $2 \leq j \leq n-i$.
 Then
 \begin{eqnarray*}
 \sum_{j=2}^{n-i} jq^{1/2}_j & \leq &  \sum_{j=J+1}^\infty j
\eta^{(j-1)/2} + \sum_{j=2}^J jq^{1/2}_j \\
 & \leq &  \sum_{j=J+1}^\infty j \eta^{(j-1)/2} + J \sqrt{ \sum_{j=2}^J
jq_j}
 \end{eqnarray*} 
 by the Cauchy-Schwarz inequality.  Setting $d_i = \sum_{j=2}^J jq_j$,
another use of Cauchy-Schwarz gives 
 \[ \sum_{i=0}^{n-2} \sqrt{d_i} \leq n^{1/2}  \sqrt{ \sum_{i=0}^{n-2} d_i}
\]
 and the result follows.

\newpage
 %\bibliographystyle{plain}
% \bibliography{../../trees/me,../../trees/networks,../../trees/alg,../../trees/misc,../../trees/trees}

\begin{thebibliography}{1}

\bibitem{me-spatial-4}
D.J. Aldous.
\newblock Which connected spatial networks on random points have linear
  route-lengths?
\newblock In preparation, 2009.

\bibitem{me116}
D.J. Aldous and W.S. Kendall.
\newblock Short-length routes in low-cost networks \emph{via} {P}oisson line
  patterns.
\newblock {\em Adv. in Appl. Probab.}, 40:1--21, 2008.

\bibitem{me-spatial-1}
D.J. Aldous and J.~Shun.
\newblock Models for connected networks over random points and a route-length
  statistic.
\newblock In preparation, 2009.  Draft available at http://www.stat.berkeley.edu/$\sim$aldous/Papers/me-spatial-1.pdf

\bibitem{MR1452554}
C.~D. Howard and C.~M. Newman.
\newblock Euclidean models of first-passage percolation.
\newblock {\em Probab. Theory Related Fields}, 108(2):153--170, 1997.

\bibitem{jaromczyk}
J.W. Jaromczyk and G.T. Toussaint.
\newblock Relative neighborhood graphs and their relatives.
\newblock {\em Proceedings of the IEEE}, 80(9):1502--1517, 1992.


\bibitem{MR905330}
H.~Kesten.
\newblock Percolation theory and first-passage percolation.
\newblock {\em Ann. Probab.}, 15(4):1231--1271, 1987.

\bibitem{kesten-FPP}
H.~Kesten.
\newblock First-passage percolation.
\newblock In {\em From Classical to Modern Probability}, number~54 in Progr.
  Probab., pages 93--143. Birkhauser, 2003.

\bibitem{MR0438477}
J.~F.~C. Kingman.
\newblock Subadditive processes.
\newblock In {\em \'{E}cole d'\'{E}t\'e de {P}robabilit\'es de {S}aint-{F}lour,
  {V}--1975}, pages 167--223. Lecture Notes in Math., Vol. 539. Springer,
  Berlin, 1976.

\bibitem{MR1166620}
M.Q. Vahidi-Asl and J.~C. Wierman.
\newblock A shape result for first-passage percolation on the {V}orono\u\i\
  tessellation and {D}elaunay triangulation.
\newblock In {\em Random graphs, {V}ol.\ 2 ({P}ozna\'n, 1989)}, Wiley-Intersci.
  Publ., pages 247--262. Wiley, New York, 1992.

\end{thebibliography}

\end{document}